\documentclass[12pt,fleqn]{article}

\usepackage{a4wide}

\usepackage{amsmath,amstext}
\usepackage{amsfonts}
\usepackage{amsthm}
\usepackage{amssymb}
\usepackage{cite}
\usepackage{mathrsfs}

\newtheorem{thm}{Theorem}[section]
\newtheorem{prop}[thm]{Proposition}
\newtheorem{lm}[thm]{Lemma}
\newtheorem{cor}[thm]{Corollary}
\newtheorem{defn}[thm]{Definition}
\DeclareMathOperator{\End}{\text{End\,}}
\DeclareMathOperator{\Id}{\text{Id}}
\numberwithin{equation}{section}

\begin{document}

\baselineskip 21pt

\parskip 7pt

\hfill  \today

\vspace{24pt}

\begin{center}


  {\large
    \textbf{
      Theta Functions Associated with the Affine Root Systems
      and \\
      the Elliptic Ruijsenaars Operators
      }
    }

  \vspace{24pt}
  Yasushi \textsc{Komori}
  \footnote[1]{E-mail: \texttt{
      komori@monet.phys.s.u-tokyo.ac.jp
      }
    }

  \vspace{8pt}

  \textsl{Department of Physics, Graduate School of Science, \\
    University of Tokyo, \\
    Hongo 7-3-1, Bunkyo, Tokyo 113, Japan.}

  \vspace{8pt}


(Received: \hspace{35mm})

\end{center}

\vspace{18pt}

\begin{center}
  \underline{ABSTRACT}
\end{center}

We study a family of mutually commutative difference operators 
associated with the affine root systems.
These operators act on the space of meromorphic functions
on the Cartan subalgebra of the affine Lie algebra.
We show that the space spanned by the characters of 
a fixed positive level
is invariant under the action of these operators.

\newpage

\section{Introduction}
In \cite{Ruij87}, 
a family of mutually commutative operators,
whose coefficients consist of theta functions,
were introduced as
a relativistic quantum many-body system, i.e.,
an elliptic difference analogue of the
Calogero-Sutherland model.
Since then, these operators have been studied extensively
from various points of view, especially 
from 
the analogy with the Macdonald operators.
The eigenvectors of the Macdonald operators are a two-parameter
extension of the Schur functions or the characters of 
finite dimensional simple Lie algebras.
Then it is natural to expect this structure in the elliptic case.
In fact, it was clarified 
in \cite{Has97,HIK98} that 
the elliptic analogues of type $A^{(1)}_l$ and $C^{(1)}_2$
have an invariant subspace
in the meromorphic functions and 
that this space is actually spanned by the characters of the
corresponding affine Lie algebra.
These facts are found through 
the studies of the intertwining vectors
between the face models and the vertex models.
Independently,
in \cite{FelPas94}, the Boltzmann weight of 
the matrix elements of Belavin's
elliptic $R$-matrix was calculated making use of this fact implicitly.

In a series of Cherednik's papers \cite{Chere92b,Chere95c,Chere97},
it has been proved out that
the double affine Hecke algebra plays an essential role 
in the Macdonald theory.
There are some algebras that are considered 
to describe the structure of 
the elliptic analogues \cite{Skl83,Chere95b,Etin94,FelVar97}.
In this paper, we employ yet another approach or the root algebra
to these operators. 
Following the previous work \cite{KomHik97c} 
where we studied nontwisted cases,
we construct a family of mutually commuting 
difference operators
associated 
with arbitrary affine root systems.
These operators are shown to act on
the vector space 
of the Weyl group invariant meromorphic functions
and, furthermore, on the space spanned by the characters
of a fixed positive level.

This paper is organized as follows:
In section 2, we prepare the notations and definitions used in this paper.
In section 3, we define the root algebras that was introduced by Cherednik
in the development of the theory of the affine Hecke algebras.
In section 4,
we demonstrate some examples of the generators of
a commutative subalgebra in the root algebras.
In section 5, we give some representations of the root algebras
with a spectral parameter,
which consist of Jacobi's theta functions and act on the
meromorphic functions on the Cartan subalgebra.
We show that when we assign a special value to the spectral parameter,
the difference operators preserve the Weyl group invariant subspace.
By construction, they form a commutative family.
In section 6, we calculate the explicit forms of these operators
at this spectral parameter and
observe that they can be regarded as an elliptic analogue of
the Macdonald operators.
In twisted cases, we have the difference and quantum version of the
systems that is recently proposed and
is dealt in terms of the Lax formalism \cite{BS98c}.
We also prove that the generators are algebraically independent
and thus the commutative subalgebra is isomorphic to
a polynomial ring.
In section 7, we show the main
theorem (Theorem \ref{thm:main}) that
 they have 
an infinite dimensional invariant subspace 
(finite rank submodule) 
of the theta functions of positive level, where
the key of the proof is due to \cite{FelPas94,KomHik96a}.
The last section is devoted to the concluding remarks.

To end this section, we present 
two elliptic difference operators
which take the simplest form among the generators respectively
in the root systems of type $A^{(1)}_{l-1}$ and $A^{(2)}_{2l}$:
\begin{align}
  & \hspace*{-1cm}
\Hat{Y}_{A^{(1)}_{l-1}}^{-\lambda_1}=\sum_{j=1}^l\prod_{j\neq k}^l
   \frac{\vartheta_1(x_j-x_k-\mu)}{\vartheta_1(x_j-x_k)}
  t_j(\kappa)
 \prod_{j=1}^l t_j(-\kappa/l), 
\label{A1_l-1}
  \displaybreak[0] \\  
  & \hspace*{-1cm}
\Hat{Y}_{A^{(2)}_{2l}}^{-\lambda_1}
  =
  \sum_{j=1}^l
    \Bigl(
      \prod_{\substack{k = 1 \\ k \neq j } }^l
      \frac{\vartheta_1(x_j-x_k-\mu)}{\vartheta_1(x_j-x_k)}
      \frac{\vartheta_1(x_j+x_k-\mu)}{\vartheta_1(x_j+x_k)}
    \Bigr) \,
    \Bigl(
      \prod_{r=0}^3
      \frac{\vartheta_{r+1}(x_j-\nu_r)}{\vartheta_{r+1}(x_j)}
      \frac{\vartheta_{r+1}(x_j+\kappa/2-\Bar{\nu}_r)}
      {\vartheta_{r+1}(x_j+\kappa/2)}
    \Bigr) 
  t_j ( \kappa )
  \notag \displaybreak[0] \\ 
  & \hspace*{-1cm}
  + 
  \sum_{j=1}^l
    \Bigl(
      \prod_{\substack{k = 1 \\ k \neq j } }^l
      \frac{\vartheta_1(-x_j-x_k-\mu)}{\vartheta_1(-x_j-x_k)}
      \frac{\vartheta_1(-x_j+x_k-\mu)}{\vartheta_1(-x_j+x_k)}
    \Bigr) \,
    \Bigl(
      \prod_{r=0}^3
      \frac{\vartheta_{r+1}(-x_j-\nu_r)}{\vartheta_{r+1}(-x_j)}
      \frac{\vartheta_{r+1}(-x_j+\kappa/2-\Bar{\nu}_r)}
      {\vartheta_{r+1}(-x_j+\kappa/2)}
    \Bigr) 
  t_j( -\kappa )
  \notag \displaybreak[0] \\
  & 
    -\sum_{p=0}^3 
    \Bigl(
      \frac {\pi}{\vartheta_1'(0)}
    \Bigr)^2
    \frac{2
      }
    {\vartheta_1(-\mu)\vartheta_1(-\kappa-\mu)}
    \Bigl(
       \prod_{r=0}^3 
       \vartheta_{r+1}(-\kappa-\nu_{\pi_p r})    
       \vartheta_{r+1}(-\Bar{\nu}_{\pi_p r})    
    \Bigr)
    \notag \displaybreak[0] \\
    & \qquad \qquad \times 
    \Bigl(
      \prod_{j=1}^N
      \frac{\vartheta_{p+1}(x_j-\kappa/2-\mu)}{\vartheta_{p+1}(x_j-\kappa/2)}
      \frac{\vartheta_{p+1}(-x_j-\kappa/2-\mu)}{\vartheta_{p+1}(-x_j-\kappa/2)}
    \Bigr) .
\label{A2_2l}
\end{align}
Here
we have realized the root systems in $\mathbb{C}^l$ in the standard way,
$\vartheta_j(x)=\vartheta_j(x;\tau)$ is the Jacobi theta function and
$t_i(\kappa)$ is a translation of the variable $x_i$ by $\kappa$.
$\pi_r$ ($r=0,1,2,3$) denotes the
permutation:
$\pi_0 = id$,
$\pi_1 = (0 1) (2 3)$,
$\pi_2 = (0 2) (1 3)$,
and
$\pi_3 = (0 3) (1 2)$.
The parameters $\kappa$, $\mu$, $\nu_r$ and $\Bar{\nu}_r$ ($r=0,1,2,3$)
are arbitrary constants.
The operator \eqref{A1_l-1} was introduced in \cite{Ruij87}
together with
the whole family of commuting difference operators,
while 
the operator \eqref{A2_2l} was 
conjectured to be a member of a commutative family
in \cite{Diej94a,Diej95a}.
The $A^{(2)}_{2l}$-type model
was referred to as $D$-type or $BC$-type in previous papers.

If we set $\kappa=l\mu/k$ in $A^{(1)}_{l-1}$ case and
$\kappa=(\nu+2\Bar{\nu}+2(l-1)\mu)/k$ in $A^{(2)}_{2l}$ case, 
where
$\nu=\sum \nu_r$ and
$\Bar{\nu}=(\sum \Bar{\nu}_r)/2$,
then these operators have an invariant subspace which consists of
the characters of level $k$ 
corresponding
to each affine Lie algebra.
When the parameters $\mu,\nu,\Bar{\nu}$ are set to be unity,
we see that $\kappa$ reduces to $h^\vee/k$ where
$h^\vee$ is the dual Coxeter number.
For the derivation of these facts
in each case, 
see \cite{Has97,KomHik97a,KomHik97b,HikKom98a,HikKom98b}.


\section{Affine Root Systems}
We give some well-known facts 
about the affine root systems and 
the affine Weyl groups~\cite{Bourbaki,Hum,Chere92a},
which are the standard tools in the theory of the affine Hecke algebras.
Some of the definitions are slightly changed and extended 
so that they include the twisted affine root systems.
Most of the notations are due to \cite{Kac}.

Let $\mathfrak{g}=\mathfrak{g}(A)$
be the affine Lie algebra 
associated with the generalized Cartan matrix $A$  
of type $X^{(r)}_N$,
$\mathfrak{h}$ its Cartan subalgebra,
$\text{dim}\,\mathfrak{h}=l+1$ the rank of $\mathfrak{g}$,
$I=\{0,\ldots,l\}$ a set of indices,
$\Pi=\{\alpha_i|i\in I\}\subset\mathfrak{h}^*$ 
the set of simple roots,
$\Pi^\vee=\{\alpha_i^\vee|i\in I\}\subset\mathfrak{h}$
the set of simple coroots,
$\Delta$ the root system,
$Q$ and $Q^\vee$ the root and coroot lattices,
$P$ and $P^\vee$ the weight and coweight lattices:
\begin{gather}
  Q=\bigoplus_{i\in I}\mathbb{Z}\,\alpha_i\subset
  P=\bigoplus_{i\in I}\mathbb{Z}\,\Lambda_i\oplus\mathbb{C}\delta 
    \subset\mathfrak{h}^*  ; \\
  Q^\vee=\bigoplus_{i\in I}\mathbb{Z}\,\alpha_i^\vee\subset
  P^\vee=\bigoplus_{i\in I}\mathbb{Z}\,\Lambda_i^\vee\oplus\mathbb{C}K
    \subset\mathfrak{h},
\end{gather}
where 
$\langle\alpha_i,\Lambda_j^\vee\rangle=\delta_{ij}$,
$\langle\Lambda_i,\alpha_j^\vee\rangle=\delta_{ij}$,
$d=\Lambda_0^\vee$.
Since the normalized invariant form is nondegenerate on
$\mathfrak{h}$,
we have an isomorphism $\nu:\mathfrak{h}\rightarrow\mathfrak{h}^*$
defined by
\begin{equation}
  \langle\nu(h),h_1\rangle=(h|h_1),\quad h,h_1\in\mathfrak{h},
\end{equation}
and the induced bilinear form $(.|.)$ on $\mathfrak{h}{}^*$.
Let $\overset{\circ}{I}=\{1,\ldots,l\}$,
$\overset{\circ}{\Pi}=\{\alpha_i|i\in\overset{\circ}{I}\}$ and 
$\overset{\circ}{\Pi}{}^\vee=\{\alpha_i^\vee|i\in\overset{\circ}{I}\}$.
Let $\overset{\circ}{\mathfrak{h}}{}^*$ be 
the subspace of $\mathfrak{h}^*$ spanned by $\overset{\circ}{\Pi}$ 
over $\mathbb{C}$.
For $\lambda\in\mathfrak{h}^*$,
denote by $\overline{\lambda}$ 
the orthogonal projection of $\lambda$ on 
$\overset{\circ}{\mathfrak{h}}{}^*$.
Let
$\overset{\circ}{Q}$ 
be the sublattice of $Q$ generated by
$\overset{\circ}{\Pi}$ and
$\overset{\circ}{P}$ 
the projection of $P$ on
$\overset{\circ}{\mathfrak{h}}{}^*$.
The dual notions 
$\overset{\circ}{\mathfrak{h}}$,
$\overline{h}$,
$\overset{\circ}{Q}{}^\vee$
and
$\overset{\circ}{P}{}^\vee$ 
are defined similarly:
\begin{gather}
  \overset{\circ}{Q}=\bigoplus_{i\in\overset{\circ}{I}}\mathbb{Z}\,\alpha_i\subset
  \overset{\circ}{P}=\bigoplus_{i\in\overset{\circ}{I}}\mathbb{Z}\,\overline{\Lambda_i}\subset
  \overset{\circ}{\mathfrak{h}}{}^*  ; \\
  \overset{\circ}{Q}{}^\vee=\bigoplus_{i\in\overset{\circ}{I}}\mathbb{Z}\,\alpha_i^\vee\subset
  \overset{\circ}{P}{}^\vee=\bigoplus_{i\in\overset{\circ}{I}}\mathbb{Z}\,\overline{\Lambda_i^\vee}\subset
  \overset{\circ}{\mathfrak{h}}.
\end{gather}
For $\alpha\in\Delta^{re}$, let $r_\alpha$ be a reflection defined by
\begin{equation}
  r_\alpha(\lambda):=\lambda-\langle\lambda,\alpha^\vee\rangle\alpha,
\qquad \lambda\in\mathfrak{h}^*.
\end{equation}
The Weyl group $\overset{\circ}{W}$ is generated 
by the fundamental reflections 
$\{r_i:=r_{\alpha_i} | i\in \overset{\circ}{I} \}$
on $\mathfrak{h}^*$
and
the affine Weyl group $W$ is generated 
by $\{r_i | i\in I \}$.
The defining relations are given by $r_i^2=id$ and the Coxeter relations: 
\begin{equation}
  (r_i\,r_j)^{m_{ij}}=id, \qquad 
  \text{for } i\neq j\in I,
\end{equation}
where 
$m_{ij}=2$ if $\alpha_i$ and $\alpha_j$
are disconnected in the Dynkin diagram $S(A)$ and
$m_{ij}=3,4,6$ if 1,2,3 lines respectively connect 
$\alpha_i$ and $\alpha_j$ in $S(A)$.
We note that there is no Coxeter relation in the affine root systems of
rank $2$.
For $\alpha\in\overset{\circ}{\mathfrak{h}}{}^*$,
we define endomorphisms $t_\alpha,t^\iota_\alpha$ 
of the vector space $\mathfrak{h}^*$ for $\kappa\in\mathbb{C}$
by (cf.~\cite{Chere95b})
\begin{gather}
t_\alpha(\lambda)
:=\lambda+\langle\lambda,K\rangle\alpha-
((\lambda|\alpha)+\frac{1}{2}|\alpha|^2
\langle\lambda,K\rangle)\delta, \\
t^\iota_\alpha(\lambda)
:=\lambda-\kappa(\lambda|\alpha)\delta.
\end{gather}
Here
$t^\iota_\alpha$ is associated with an endomorphism
of $\mathfrak{h}{}^*$,
$\iota(\lambda):=\lambda+(\kappa-1)\langle\lambda,d\rangle\delta$
as follows:
\begin{equation}
\iota \circ t_\alpha (\overset{\circ}{\lambda}+m\delta)= 
t^\iota_\alpha (\overset{\circ}{\lambda}+m\delta).
\end{equation}
Let $a_i$ and $a_i^\vee$
be the labels of the Dynkin diagram from Table Aff
in \cite{Kac}.
Note that $a_0=2$ if $A$ is of type $A^{(2)}_{2l}$ and $a_0=1$ otherwise.
Let $\theta:=\delta-a_0\alpha_0 \in\overset{\circ}{\Delta}_+$,
$M:=\nu(\mathbb{Z}(\overset{\circ}{W}\cdot\theta^\vee))
\subset\overset{\circ}{\mathfrak{h}}{}^*$, 
$T_M$ the corresponding group of translations of $M$. 
Then
\begin{prop}
The group $W$ is the semidirect product
$W=\overset{\circ}{W}\ltimes T_M$.  
\end{prop}
For $\alpha\in\Delta^{re}$,
let $\gamma_\alpha:=r$ if $\alpha\in\Delta_l$ 
and    
$\gamma_\alpha:=1$ otherwise.
Then the real roots are written as
\begin{equation}
  \Delta^{re}=
  \begin{cases}
  \{\alpha+n\gamma_\alpha\delta|
  \alpha\in\overset{\circ}{\Delta},n\in\mathbb{Z}\},
  &\text{if $A$ is not of type $A^{(2)}_{2l}$};\\[2mm]
  \begin{split}
  &\{\alpha+n\gamma_\alpha\delta|
  \alpha\in\overset{\circ}{\Delta},n\in\mathbb{Z}\}\cup \\
  &\quad\{\frac 12(\alpha+(2n-1)\delta)|
  \alpha\in\overset{\circ}{\Delta}_l,n\in\mathbb{Z}\},
  \end{split}
  \qquad&\text{if $A$ is of type $A^{(2)}_{2l}$}.  
  \end{cases}
\end{equation}
Let $\widehat{M}:=\{\lambda\in\overset{\circ}{\mathfrak{h}}{}^*|
\alpha\in\Delta^{re},(\alpha|\lambda)\in\gamma_\alpha\mathbb{Z}\}$.
Then we see that 
$\widehat{M}\subset\overset{\circ}{P}$ and
$T_{\widehat{M}}$ is normalized by
$\overset{\circ}{W}$.
\begin{defn}
The extended affine Weyl group $\widehat{W}$ is the semidirect product
$\widehat{W} := \overset{\circ}{W}\ltimes T_{\widehat{M}}$.
\end{defn}
The lattice $\widehat{M}$ is defined so that 
the extended affine Weyl group acts on $\Delta$.
Here are 
the explicit description of $\widehat{M}$ and
its canonical basis $\{\lambda_i|i\in\overset{\circ}{I}\}$:
\begin{equation}
  \widehat{M}=
  \begin{cases}
   \nu(\overset{\circ}{P}{}^\vee),\qquad &\text{if $r=1$}; \\
   \overset{\circ}{P}, &\text{otherwise,}
  \end{cases}
\qquad
  \lambda_i=
  \begin{cases}
    \nu(\overline{\Lambda_i^\vee}),\qquad &\text{if $r=1$}; \\[1.5mm]
    \overline{\Lambda_i}, &\text{otherwise}.
  \end{cases}
\end{equation}
We also use 
$\widehat{M}_-:=\oplus_{i\in\overset{\circ}{I}}\mathbb{Z}_{\leq 0}\lambda_i$.
The action of $\widehat{W}$ is naturally induced on $\mathfrak{h}$
via the form $\langle\cdot,\cdot\rangle$.

Let $\Omega$ be the subgroup of $\widehat{W}$
which stabilizes the affine Weyl chamber $C$. 
\begin{prop}
The subgroup $\Omega$ is isomorphic 
to $\widehat{W}/W \simeq T_{\widehat{M}}/T_{M}$ thus Abelian.
The extended affine Weyl group $\widehat{W}$
is isomorphic to the semi-direct product $W\rtimes\Omega$.
\end{prop}
\begin{defn}
\begin{enumerate}
\item The length $\ell(w)$ of $w\in W$ is the length $\ell$ of 
the reduced decomposition:
\begin{gather}
  w=r_{i_1}\ldots r_{i_\ell},\qquad \text{for}\qquad i_k\in I, \\
  \ell(id)=0.
\end{gather}
\item The length $\ell(\Hat{w})$ of $\hat{w}\in \widehat{W}$ is the
number of the positive roots made negative by $\Hat{w}^{-1}$:
\begin{gather}
  \ell(\Hat{w}):=|\Delta_{\Hat{w}}|, \\
  \Delta_{\Hat{w}}:=\{\alpha\in \Delta_+\cap -\Hat{w}\Delta_+\},
\end{gather}
which is equivalent to
the definition $\ell(w)$ for $w\in W$.
The reduced decomposition of 
$\Hat{w}\in\widehat{W}$ is 
$\Hat{w}=w\omega=r_{i_1}\ldots r_{i_\ell}\omega$,
where
$\omega\in\Omega$ and $\ell=\ell(\Hat{w})=\ell(w)$.
\end{enumerate}
\end{defn}
The set $\Delta_{\Hat{w}}$ is  
explicitly described as $\Delta_{\Hat{w}}=\{
\alpha^1=\alpha_{i_1},\alpha^2=r_{i_1}(\alpha_{i_2}),\ldots,
\alpha^{\ell}=w r_{i_\ell}(\alpha_{i_\ell})\}$.
By definition, $\Delta_{\Hat{w}}$ is independent of
reduced expressions.
One sees that $\Omega=\{\omega \in \widehat{W}, \ell(\omega)=0\}$.
\begin{defn}
A weight $\lambda\in\widehat{M}$ is said to be minuscule if 
$\Delta_{t_{-\lambda}}\subset\overset{\circ}{\Delta}_+$.
\end{defn}

We use the following useful formulas,
which can be easily derived from the definitions above:
\begin{subequations}
\label{property}
\begin{gather}
  \Delta_{t_{\lambda_-}}=
  \begin{cases}
        \{\alpha-n\gamma_\alpha\delta |
          \alpha\in\overset{\circ}{\Delta}_+,   
        0 \geq n > \dfrac {1}{\gamma_\alpha}(\lambda_-|\alpha)  \} , \qquad
  &\text{if $A$ is not of type $A^{(2)}_{2l}$};\\[5mm]
  \begin{split}
  &      \{\alpha-n\gamma_\alpha\delta |
          \alpha\in\overset{\circ}{\Delta}_{+},   
        0 \geq n > \frac {1}{\gamma_\alpha}(\lambda_-|\alpha)  \}\cup \\
  &\quad  \{\frac 12(\alpha-(2n-1)\delta) |
          \alpha\in(\overset{\circ}{\Delta}_{+})_l,   
        0 \geq n > \frac 12(\lambda_-|\alpha)  \},
  \end{split}
  &\text{if $A$ is of type $A^{(2)}_{2l}$},
  \end{cases}
\displaybreak[0]
\\
  \ell(t_{\lambda_-})=
  \begin{cases}
   \sum_{\alpha\in\overset{\circ}{\Delta}_+}
   \Bigl|\frac{1}{\gamma_\alpha}(\alpha|\lambda_-)\Bigr|, \qquad 
   &\text{if $A$ is not of type $A^{(2)}_{2l}$};\\[2mm]
   \sum_{\alpha\in\overset{\circ}{\Delta}_+}
   \Bigl|(\alpha|\lambda_-)\Bigr|, \qquad 
   &\text{if $A$ is of type $A^{(2)}_{2l}$}, 
  \end{cases}
\displaybreak[0]
\\
  \ell(r_j\,t_{-\lambda_i})=\ell(t_{-\lambda_i})+1,
\displaybreak[0]
\\
  \ell(r_i\,t_{-\lambda_i})=\ell(t_{-\lambda_i})-1,
\displaybreak[0]
\\
  \ell(t_{\lambda_-}\,w)=\ell(t_{\lambda_-})+\ell(w),
\displaybreak[0]
\\
  \ell(t_{\lambda_-+\lambda'_-})=\ell(t_{\lambda_-})+\ell(t_{\lambda'_-}),
\end{gather}
\end{subequations}
where $i\neq j\in\overset{\circ}{I}$, $\lambda_-,\lambda'_-\in\widehat{M}_-$, 
$w \in \overset{\circ}{W}$.

\section{Root Algebras}
We shall define the root algebras after 
Cherednik \cite{Chere92a}.
Let $\mathcal{T}$ be the tensor algebra over $\mathbb{C}$
generated
by independent variables
$\{R_\alpha | \alpha \in \Delta^{re}\}$.
Then the action of $\Hat{w}\in\widehat{W}$ on $\Delta^{re}$
induces an action on $\mathcal{T}$ 
by ${}^{\Hat{w}}:R_{\alpha} \mapsto R_{\Hat{w}(\alpha)}$.
\begin{defn}
Let $\mathcal{I}$ be the ideal in $\mathcal{T}$ 
which is generated by
all the elements of the form for
$i\neq j\in I$,
and $\Hat{w}\in \widehat{W}$ :
\begin{equation}
\label{rel_Yang}
  {}^{\Hat{w}}(
\underbrace{R_{\alpha_i}\otimes R_{r_i \alpha_j}\otimes R_{r_i r_j \alpha_i}\otimes\cdots}_{m_{ij}\text{ factors}}
)
- {}^{\Hat{w}}(
\underbrace{R_{\alpha_j}\otimes R_{r_j \alpha_i}\otimes R_{r_j r_i \alpha_j}\otimes \cdots}_{m_{ij}\text{ factors}}
).
\end{equation}
The root algebra $\widetilde{\mathcal{R}}$ is $\mathcal{T}/\mathcal{I}$.
$\{R_\alpha | \alpha \in \Delta^{re}\}$ are called the $R$-matrices. 
\end{defn}
Because of the $\widehat{W}$-invariance of $\mathcal{I}$,
the action of $\widehat{W}$ is induced on $\widetilde{\mathcal{R}}$.
For simplicity,
we write products in $\widetilde{\mathcal{R}}$ in the usual way 
for associative algebras.
\begin{thm}
\begin{enumerate}
\item There exists a unique set 
$\{R_{\Hat{w}} | \Hat{w} \in \widehat{W}\}\subset\widetilde{\mathcal{R}}$ 
satisfying the relations:
\begin{equation}
  R_{v\,w}=R_v\,^{v}\!R_w,\qquad R_{r_i}=R_{\alpha_i}
\quad(i\in I), 
\qquad R_{\omega}=1,
\end{equation}
where $\omega\in\Omega$, $v,w \in \widehat{W}$ and 
$\ell(v\,w)=\ell(v)+\ell(w)$. 
\item
We have the $R$-matrix
for $\Hat{w}\in\widehat{W}$ and its arbitrary reduced decomposition 
$\Hat{w}= w \omega =r_{i_1} \ldots r_{i_\ell}\omega $ as
\begin{equation}
  R_{\Hat{w}}=R_{\alpha^1}\ldots R_{\alpha^\ell},\qquad
\alpha^{1}=\alpha_{i_1},\quad
\alpha^{2}=r_{i_1}(\alpha_{i_2}),\quad
\ldots,\quad
\alpha^{\ell}=w r_{i_\ell}(\alpha_{i_\ell})\in \Delta_{\Hat{w}}.
\end{equation}
\end{enumerate}
\end{thm}
Instead of the original root algebra,
we use the following extension, where 
$\widetilde{\mathcal{R}}$ is combined with 
the translation group $T_{\widehat{M}}$:
\begin{defn}
\label{R_sp_t}
$\mathcal{R}:=\widetilde{\mathcal{R}}\rtimes T_{\widehat{M}}$:
\begin{equation}
  (R\,t_\lambda)(R'\,t_\mu)
 =R\,({}^{t_\lambda}\!R')\,t_{\lambda+\mu},
\end{equation}
where $R,R'\in\widetilde{\mathcal{R}}$ and $\lambda,\mu\in\widehat{M}$.
\end{defn}
We see that $\mathcal{R}$ is generated by
$\{t_{\lambda_i},R_{\alpha}|
i\in\overset{\circ}{I},\alpha\in\overset{\circ}{\Delta}\}$ 
if $A$ is not of type $A^{(2)}_{2l}$ and  
$\{t_{\lambda_i},R_{\alpha}|
i\in\overset{\circ}{I},\alpha\in\overset{\circ}{\Delta},
2\alpha-\delta\in\overset{\circ}{\Delta}\}$ 
if $A$ is of type $A^{(2)}_{2l}$. 

\begin{thm}
The subalgebra $\mathcal{S}\subset\mathcal{R}$ generated by 
$\{Y^\lambda:=R_{t_\lambda}\,t_\lambda |
\lambda\in\widehat{M}_-\}$ 
forms a commutative algebra and is generated by
$\{Y^{-\lambda_i} | i\in\overset{\circ}{I}\}$.
\end{thm}
\begin{proof}
It is straightforward by the formulas \eqref{property} and Definition
\ref{R_sp_t}.
\end{proof}
 
\section{Affine Root Systems of Rank $3$}
There are six types of affine root systems of rank $3$.
We denote
$\alpha= \alpha_1$ and $\beta = \alpha_2$ 
where $|\alpha_1|\geq|\alpha_2|$,
and $\lambda=\lambda_1$ and $\mu=\lambda_2$, respectively.
We have the following systems:
\begin{subequations}
\begin{align}
\text{$A^{(1)}_2$-type} \qquad &
Y^{-\lambda}=R_{\alpha}\,R_{\alpha+\beta}\,t_{-\lambda} ,\\
&Y^{-\mu}=R_{\beta}\,R_{\alpha+\beta}\,t_{-\mu} ,
\end{align}
\end{subequations}
\begin{subequations}
\begin{align}
\text{$C^{(1)}_2$-type} \qquad &
Y^{-\lambda}=R_{\alpha}\,R_{\alpha+\beta}
\,R_{\alpha+2\beta}\,t_{-\lambda} ,\\
&Y^{-\mu}=R_{\beta}\,R_{\alpha+2\beta}\,R_{\alpha+\beta}
\,R_{\alpha+2\beta+\delta}\,t_{-\mu} ,
\end{align}
\end{subequations}
\begin{subequations}
\begin{align}
\text{$G^{(1)}_2$-type} \qquad &
Y^{-\lambda}=R_{\alpha}\,R_{\alpha+\beta}
\,R_{2\alpha+3\beta}\,R_{\alpha+2\beta}
\,R_{\alpha+3\beta}\,R_{2\alpha+3\beta+\delta}
\,t_{-\lambda} ,\\[2mm]
\begin{split}
& Y^{-\mu}=R_{\beta}\,R_{\alpha+3\beta}\,R_{\alpha+2\beta}
\,R_{2\alpha+3\beta}\,R_{\alpha+\beta}
\\[-1mm]
& \qquad \quad R_{\alpha+3\beta+\delta}\,R_{2\alpha+3\beta+\delta}
\,R_{\alpha+2\beta+\delta}\,R_{\alpha+3\beta+2\delta}
\,R_{2\alpha+3\beta+2\delta}
\,t_{-\mu} ,
\end{split}
\end{align}
\end{subequations}
\begin{subequations}
\begin{align}
\text{$A^{(2)}_4$-type} \qquad 
&Y^{-\lambda}=R_{\alpha}\,R_{\alpha+\beta}\,R_{\alpha+2\beta}
\,R_{\frac 12\alpha+\frac 12\delta}\,R_{\alpha+\beta+\delta}
\,R_{\frac 12\alpha+\beta+\frac 12\delta}\,t_{-\lambda} ,\\
&Y^{-\mu}=R_{\beta}\,R_{\alpha+2\beta}\,R_{\alpha+\beta}
\,R_{\frac 12\alpha+\beta+\frac 12\delta}\,t_{-\mu} ,
\end{align}
\end{subequations}
\begin{subequations}
\begin{align}
\text{$D^{(2)}_3$-type} \qquad 
&Y^{-\lambda}=R_{\alpha}\,R_{\alpha+\beta}\,R_{\alpha+2\beta}
\,R_{\alpha+\beta+\delta}\,t_{-\lambda} ,\\
&Y^{-\mu}=R_{\beta}\,R_{\alpha+2\beta}\,R_{\alpha+\beta}\,t_{-\mu} ,
\end{align}
\end{subequations}
\begin{subequations}
\begin{align}
\text{$D^{(3)}_4$-type} \qquad 
\begin{split}
&Y^{-\lambda}=R_{\alpha}\,R_{\alpha+\beta}\,R_{2\alpha+3\beta}\,
R_{\alpha+2\beta}\\[-1mm]
& \qquad \quad R_{\alpha+3\beta}\,R_{\alpha+\beta+\delta}\,
R_{\alpha+2\beta+\delta}\,R_{2\alpha+3\beta+3\delta}\,
R_{\alpha+\beta+2\delta}\,R_{\alpha+2\beta+2\delta}\,t_{-\lambda} ,
\end{split}
\\[2mm]
&Y^{-\mu}=R_{\alpha}\,R_{\alpha+3\beta}\,R_{\alpha+2\beta}\,
R_{2\alpha+3\beta}\,R_{\alpha+\beta}\,R_{\alpha+2\beta+\delta}\,t_{-\mu} .
\end{align}
\end{subequations}

\section{Elliptic $R$-matrices}
For $\alpha\in\Delta^{re}$, 
let $\mu_\alpha\in\mathbb{C}$ be $\widehat{W}$-invariant constants:
$\mu_{\hat{w}(\alpha)}=\mu_\alpha$ for $\hat{w}\in\widehat{W}$.
Let $Y:=\{h\in\mathfrak{h}\,|\ \text{Re}\langle \delta,h \rangle>0 \}$ and
let $\mathcal{M}$ be the set of
meromorphic functions on $Y$.
We define
an action of $w=\overset{\circ}{w}t_\lambda
\in\widehat{W}$ 
on $\mathcal{M}$
as
$(w\,f)(h):=f(t^\iota_{-\lambda}\overset{\circ}{w}{}^{-1} (h))$. 

Fix $\kappa\in\mathbb{C}$ 
and $\xi\in\overset{\circ}{\mathfrak{h}}{}^*$.
We define $\Hat{R}_\alpha\in\End_\mathbb{C}\mathcal{M}$ for $\alpha\in\Delta^{re}$ by
\begin{equation}
\label{R_matrix}
\Hat{R}_\alpha
 :=
  H_\alpha(\mu_\alpha) 
  -
  H_\alpha(\langle\xi,\alpha^\vee\rangle) \, r_{\alpha} 
  ,
\end{equation}
with the following function
(see the definitions in Appendix):
\begin{equation}
\label{Solution_R_1}
  H_\alpha(\nu) 
  := 
\frac{
\vartheta^1(-\gamma_\alpha\mu_\alpha\delta;\gamma_\alpha) }
{ \vartheta^{1\prime}(0;\gamma_\alpha) }
  \sigma_{\gamma_\alpha\nu}(\iota(\alpha);\gamma_\alpha).
\end{equation}
\begin{thm}
\label{Thm_R_1}
The map $\pi:R_\alpha \mapsto \Hat{R}_\alpha, t_\lambda \mapsto t_\lambda$
induces a homomorphism from $\mathcal{R}$ to $\End_\mathbb{C}\mathcal{M}$.
These $R$-matrices satisfy the
unitarity 
\begin{equation}
\Hat{R}_\alpha\,\Hat{R}_{-\alpha}=
\Bigl(
\frac{
\vartheta^1(-\gamma_\alpha\mu_\alpha\delta;\gamma_\alpha) }
{ \vartheta^{1\prime}(0;\gamma_\alpha) }
\Bigr)^2
\bigl( 
\wp^0 (\gamma_\alpha\mu_\alpha\delta ;\gamma_\alpha) -
\wp^0 (\gamma_\alpha\langle\xi,\alpha^\vee\rangle\delta ;\gamma_\alpha) 
\bigr)\Id_\mathcal{M}.
\end{equation}
\end{thm}

Besides the above representation,
we have more general forms that
depend on the relation among $Q,Q^\vee,M$.
For $\alpha\in\Delta^{re}$, 
let 
\begin{equation}
\label{cond_theta}
N_\alpha
:=
\Biggl\{
\phi_\alpha^j:=(m_\alpha,n_\alpha)\in\mathbb{R}_{>0}^2 \Biggm| 
\begin{split}
&\overline{\alpha^\vee}\in m_\alpha\overset{\circ}{Q}{}^\vee,\qquad
m_\alpha\langle \alpha,\overset{\circ}{Q}{}^\vee\rangle\subset\mathbb{Z}, \\
&n_\alpha\gamma_\alpha\overline{\alpha^\vee}\in m_\alpha M,
m_\alpha(M|\alpha)\subset n_\alpha\gamma_\alpha\mathbb{Z} 
\end{split}
\Biggr\}.
\end{equation}
This condition is required when the root algebra acts on the 
vector space spanned by theta functions
(Proposition \ref{prop:R_end_Th})
and is an elliptic analogue in 
the representation of
the Hecke algebras \cite{Lus89}.

We enumerate
the set $N_\alpha$ 
 as
\begin{center}
\begin{tabular}{|l||c|c|c|c|}
\hline
 & $\phi_\alpha^1$ & $\phi_\alpha^2$ & $\phi_\alpha^3$ & $\phi_\alpha^4$ \\
\hline
\hline
$A$ is of type $C^{(1)}_l$ and $\alpha$ is long &
 (1,1) & (1,2) & (1/2,1) & (1/2,1/2) \\ 
$A$ is of type $A^{(2)}_{2l-1}$ and $\alpha$ is long &
 (1,1) & & & (1/2,1/2) \\
$A$ is of type $D^{(2)}_{l+1}$ and $\alpha$ is short &
 (1,1) & (1,2) & & \\
$A$ is of type $A^{(2)}_{2l}$ and $\alpha$ is short &
 (2,1) & (2,2) & (1,1) & (1,1/2) \\
$A$ is of type $A^{(2)}_{2l}$ and $\alpha$ is long &
 (1,1/2) & (1,1) & (1/2,1/2) & (1/2,1/4) \\
otherwise &
 (1,1) & & & \\
\hline
\end{tabular}
\end{center}
Here
we have numbered the elements of $N_\alpha$ for later convenience.
Let $\zeta^j_\alpha\in\mathbb{C}$ for $1\leq j\leq 4$ 
$\widehat{W}$-invariant constants.
If $\phi_\alpha^j\not\in N_\alpha$, set $\zeta^j_\alpha=0$.
In place of \eqref{Solution_R_1},
we define
\begin{equation}
\label{Solution_R_2}
  H_\alpha(\nu)
  := 
\sum_{\phi_\alpha^j=(m_\alpha,n_\alpha)\in N_\alpha} \zeta^j_\alpha
\frac{
\vartheta^1(-n_\alpha\gamma_\alpha\mu_\alpha\delta/m_\alpha;n_\alpha\gamma_\alpha) }
{ \vartheta^{1\prime}(0;n_\alpha\gamma_\alpha) }
  \sigma_{n_\alpha\gamma_\alpha\nu/m_\alpha}(\iota(m_\alpha\alpha);n_\alpha\gamma_\alpha) .
\end{equation}
Then we have a more general representation of $\mathcal{R}$ including 
Theorem \ref{Thm_R_1}.
\begin{thm}
\label{Thm_R_2}
The map $\pi$ in Theorem \ref{Thm_R_1} with \eqref{Solution_R_2}
induces a homomorphism from $\mathcal{R}$ to $\End_\mathbb{C}\mathcal{M}$.
These $R$-matrices satisfy the
unitarity
\begin{equation}
\label{unitarity}
\Hat{R}_\alpha\,\Hat{R}_{-\alpha}=
 u_\alpha(\delta)\Id_\mathcal{M},
\end{equation}
where $u_\alpha(\delta)$ depends only on $\delta$ and vanishes 
if $\langle \xi,\alpha^\vee \rangle=\pm\mu_\alpha$.

More precisely, we have 
$u_\alpha(\delta)=((p_1.\zeta_a)^2,(p_2.\zeta_a)^2,(p_3.\zeta_a)^2,(p_4.\zeta_a)^2).S.d_\alpha$,
\begin{gather}
    S = \frac 14
    \begin{pmatrix}
      1 & 0 & 0 & 0 \\
     -1 & 0 & 0 & 1 \\
     -1 & 4 & 0 & 0 \\
      1 &-4 & 4 &-1
    \end{pmatrix}
    ,
    \qquad \quad
    \zeta_\alpha
    =
    \begin{pmatrix}
      \Tilde{\zeta}^1_\alpha \\
      \Tilde{\zeta}^2_\alpha \\
      \Tilde{\zeta}^3_\alpha \\
      \Tilde{\zeta}^4_\alpha 
    \end{pmatrix} 
    ,
    \qquad \quad
    d_\alpha
    =
    \begin{pmatrix}
      d^1_\alpha \\
      d^2_\alpha \\
      d^3_\alpha \\
      d^4_\alpha 
    \end{pmatrix} ,
    \notag\\[2mm]
    p_1=(2,1,1,2), \quad
    p_2=(0,0,1,2), \quad
    p_3=(0,1,1,0), \quad
    p_4=(0,0,1,0), \notag\\
    \Tilde{\zeta}^j_\alpha
    =
    \zeta^j_\alpha 
    \frac 
    { \vartheta^1(-n_\alpha\gamma_\alpha\mu_\alpha\delta/m_\alpha;n_\alpha\gamma_\alpha) }
    { \vartheta^{1\prime}(0;n_\alpha\gamma_\alpha) }
      ,\notag \\
    d^j_\alpha
    =
\wp^0 (n_\alpha\gamma_\alpha\mu_\alpha\delta/m_\alpha;n_\alpha\gamma_\alpha) -
\wp^0 (n_\alpha\gamma_\alpha\langle\xi,\alpha^\vee\rangle\delta/m_\alpha;n_\alpha\gamma_\alpha) 
    \notag
\end{gather}
\end{thm}
\begin{proof}
We can verify the relations \eqref{rel_Yang} case-by-case,
by a direct substitution of \eqref{Solution_R_2};
for details, see \cite{ShiUen92,KomHik96a,Kom99b}. 
\end{proof}
We employ these operators even for
the affine root systems of rank $2$,
though they do not have any Coxeter relations.

We shall clarify some properties of the operators
$\Hat{Y}^\lambda=\pi(Y^\lambda)$.

\begin{lm}
\label{commute}
The $R$-matrices $\Hat{R}$ satisfy 
the following relations:
\begin{gather}
\label{vanish_1}
  \Hat{R}_{-\alpha_j}\Hat{R}_{t_{-\lambda_i}}
  ={}^{r_j}\Hat{R}_{t_{-\lambda_i}}\Hat{R}_{-\alpha_j},
  \qquad\text{for $j\neq i$}; \\
\label{vanish_2}
  \Hat{R}_{t_{-\lambda_i}}=\Hat{R}_{\alpha_i} \mathscr{R},
\end{gather}
where $\mathscr{R}$ is a product of some $R$-matrices.
\end{lm}
\begin{proof}
Combining \eqref{property},  the unitarity \eqref{unitarity},
and an equality
$r_j\,t_{-\lambda_i}=t_{-\lambda_i}\,r_j$,
we obtain \eqref{vanish_1}
for generic $\xi$ and 
thus for all $\xi\in\overset{\circ}{\mathfrak{h}}{}^*$.
The form \eqref{vanish_2} is due to the fact that
$\ell(r_i\,t_{-\lambda_i})=\ell(t_{-\lambda_i})-1$
implies the exchange condition \cite{Bourbaki},
$t_{-\lambda_i}=r_{i_1}\ldots r_{i_\ell}\omega
=r_i r_{i_1}\ldots r_{i_{m-1}}r_{i_{m+1}}\ldots r_{i_\ell}\omega$ for some $m$.
\end{proof}

If the parameter $\xi$ 
satisfies $\langle\xi,\alpha_i^\vee\rangle=-\mu_{-\alpha_i}$,
then
the $R$-matrix $\Hat{R}_{-\alpha_i}$ reduces to the form,
$\Hat{R}_{-\alpha_i}=2\,H_{-\alpha}(\mu_{-\alpha}) P_i^{(-)}$
where $P_i^{(-)}$ is
the antisymmetric projection $\frac{1}{2}(1-r_i)$. 
Let
\begin{equation}
\overset{\circ}{\rho}_\mu
:=
\sum_{i\in\overset{\circ}{I}}\mu_{\alpha_i}\overline{\Lambda_i}
=\frac 12\sum_{\alpha\in\overset{\circ}{\Delta}_+}\mu_\alpha\alpha.  
\end{equation}

From these properties,
we have the following theorem:
\begin{thm}
\label{symmetric}
Let $\mathcal{V}:=\mathcal{M}^{\overset{\circ}{W}}$, the 
$\overset{\circ}{W}$-invariant subspace of $\mathcal{M}$
and let $\xi=-\overset{\circ}{\rho}_\mu$. 
Then $\Hat{Y}^\lambda\in\End_\mathbb{C}\mathcal{V}$.
\end{thm}
\begin{proof}
It is sufficient to check it for the generators $\Hat{Y}^{-\lambda_i}$.
By Lemma \ref{commute}, we see that
$\Hat{R}_{-\alpha_j}\Hat{Y}^{-\lambda_i}|_{\mathcal{V}}=0$,
for $j\neq i$ by \eqref{vanish_1},
and for $j=i$ by \eqref{vanish_2} noting that the unitarity
\eqref{unitarity} vanishes.
Hence
$\Hat{Y}^{-\lambda_i}|_{\mathcal{V}}
=r_j\,\Hat{Y}^{-\lambda_i}|_{\mathcal{V}}
$ 
for all $j\in\overset{\circ}{I}$.
\end{proof}

The symbol $Y^\lambda$ is adopted since in a certain limit,
it reduces to
the same one up to a constant factor
as in the affine Hecke algebras, 
where $Y^\lambda$ is defined for all $\lambda\in\widehat{M}$. 
We remark that 
$\hat{Y}^\lambda$ has the inverse 
in $\End_\mathbb{C}\mathcal{M}$ for generic 
$\xi\in\overset{\circ}{\mathfrak{h}}{}^*$,
but loses its inverse 
when $\xi=-\overset{\circ}{\rho_\mu}$.

\section{Elliptic Difference Operators}

In this section, 
we calculate the explicit forms of the operators 
$\Hat{Y}^{\lambda}$
for some $\lambda$
on the space $\mathcal{V}$.
Throughout this section,
we fix $\xi=-\overset{\circ}{\rho}_\mu$.
\begin{thm} 
\label{form_min}
Let $(-\lambda)$ be minuscule.
Then we have
\begin{equation}
\Hat{Y}^{\lambda}|_\mathcal{V}=\frac{1}{\bigl|\overset{\circ}{W}_\lambda\bigr|}
\sum_{w\in\overset{\circ}{W}}w\Biggl(
\prod_{
\substack{
  \alpha\in \overset{\circ}{\Delta}_+ \\ 
  (\lambda|\alpha)=-\gamma_\alpha}}
H_\alpha(\mu_\alpha)\,t_{\lambda}
\Biggr)\Biggr|_\mathcal{V},
\label{minuscule_op}
\end{equation}
where $\overset{\circ}{W}_{\lambda}$ is the stabilizer of $\lambda$ 
in $\overset{\circ}{W}$.
\end{thm}
\begin{proof}
First notice that
$R_{t_{\lambda}}$ consists of nonaffine $R$-matrices,
$R_\alpha$ for $\alpha\in\overset{\circ}{\Delta}_+$,
because $(-\lambda)$ is minuscule.
Substituting the $R$-matrices (\ref{R_matrix}) into $\Hat{Y}^\lambda$
and expanding them, we see that every term includes
a translation operator of the form
$t_{w(\lambda)}\,w$, where 
\begin{gather}
w=r_{\alpha\{p\}}\ldots r_{\alpha\{1\}}\in \overset{\circ}{W},
\\
\alpha\{q\}:=\alpha^{m_q}, \qquad
1\leq m_p < m_{p-1} < \ldots < m_2 < m_1 \leq \ell(t_{\lambda}).   
\end{gather}
Let us show that $w(\lambda)=\lambda$ implies $w=id$.
Suppose $w(\lambda)=\lambda$ and $w\neq id$, 
then we have $\ell(w\,t_{\lambda})=\ell(t_{\lambda}\,w)$. 
From \eqref{property}, 
$\ell(t_{\lambda}\,w)=\ell(t_{\lambda})+\ell(w) > \ell(t_{\lambda})$
 while $\ell(w\,t_{\lambda})<\ell(t_{\lambda})$ by a direct calculation, 
which leads to a contradiction.
This implies that the term including $t_{\lambda}\,w$, 
$w\in\overset{\circ}{W}$ appears if and only if $w=id$.
The coefficient of this term can be easily calculated,
 \begin{equation}
\prod_{\substack{\alpha\in\overset{\circ}{\Delta}_+ \\ 
(\lambda|\alpha)=-\gamma_\alpha}}
H_\alpha(\mu_\alpha).
 \end{equation}
The $\overset{\circ}{W}$-invariance 
of the operator $\Hat{Y}^{\lambda}$ yields the form
\eqref{minuscule_op}.
\end{proof}

It is worth
noting that 
as in the trigonometric case \cite{Mac95},
we can rewrite
$\Hat{Y}^{\lambda}$ 
in a simply laced root system
as follows:
\begin{equation}
\Hat{Y}^{\lambda}|_\mathcal{V}=\frac{1}{\bigl|\overset{\circ}{W}_\lambda\bigr|}
\sum_{w\in\overset{\circ}{W}}
\frac{ (t_{-\mu w\lambda/\kappa} A_\rho) }{ A_\rho }
\,t_{w\lambda}\Biggr|_\mathcal{V},
\end{equation}
where we have set $\mu=\mu_\alpha$ and $\zeta^1_\alpha=1$.

In general,
it is complicated and difficult to compute the explicit forms of the
operators when $\lambda$ is not minuscule.
It is the case even in the framework of
the affine Hecke algebras.
There is no minuscule weight available in
the root systems of type $E^{(1)}_8, F^{(1)}_4, G^{(1)}_2, A^{(2)}_{2l},
E^{(2)}_6$ and $D^{(3)}_4$.
However, every root system possesses a ``quasi-minuscule'' weight
$\nu(\theta^\vee)$ in the sense of the following properties:
\begin{lm}
\begin{enumerate}
\label{property_theta}
\item 
$\Delta_{t_{-\nu(\theta^\vee)}}=
\Delta_{r_\theta}\cup\{a_0^{-1}(\delta+\theta)\}$.
\item
$(\nu(\theta^\vee)|\alpha)=0\text{ or }\gamma_\alpha$
for $\alpha\in\overset{\circ}{\Delta}_+$, 
$\alpha\neq\theta$, and $(\nu(\theta^\vee)|\theta)=2$.
\item
$\nu(\theta^\vee)=\lambda_i$ 
where $\alpha_i$ is the unique vertex connected to $\alpha_0$
if $A$ is not of type $A^{(1)}_l$.
\end{enumerate}
\end{lm}
\begin{proof}
We see that
$r_\theta \alpha_0=r_\theta \bigl( a_0^{-1}(\delta-\theta) \bigr)
=a_0^{-1}(\delta+\theta)\in\Delta^{re}_+$,
which implies the first statement due to the expression
$r_\theta\,r_0=t_{-\nu(\theta^\vee)}$.
The second statement is immediate from the first and \eqref{property}.
Since
$\langle\nu(\theta^\vee),\alpha_i^\vee\rangle
=\langle a_0^{-1}\theta,\alpha_i^\vee\rangle
=\langle a_0^{-1}\delta-\alpha_0,\alpha_i^\vee\rangle
=-\langle\alpha_0,\alpha_i^\vee\rangle$, 
we have
$\nu(\theta^\vee)
=-\sum_{i\in\overset{\circ}{I}}\langle\alpha_0,\alpha_i^\vee\rangle
\overline{\Lambda_i}$.
Then the last statement follows from the tables in \cite{Bourbaki,Kac}.
\end{proof}
Since in the root system of type $A^{(1)}_l$,
every $\lambda_i$ is minuscule, we have
the explicit form of $\Hat{Y}^{-\nu(\theta^\vee)}$ by 
Theorem \ref{form_min}.
So we concentrate on the other root systems.
Fix $i$ as in Lemma \ref{property_theta}. 

By the expression $t_{-\nu(\theta^\vee)}=r_\theta\,r_0$, 
we have
$Y^{-\nu(\theta^\vee)}=
R_{r_\theta}\,R_{a_0^{-1}(\theta+\delta)}\,t_{-\nu(\theta^\vee)}=
R_{r_\theta}\,t_{-\nu(\theta^\vee)}\,R_{-\alpha_0}$.
For the operator $\Hat{Y}^{-\nu(\theta^\vee)}$, an analogous statement 
to Lemma \ref{commute} holds.
\begin{lm}
The $R$-matrices $\Hat{R}$ satisfy 
the following relations:
\begin{gather}
\label{vanish_1q}
  \Hat{R}_{-\alpha_j}\Hat{R}_{r_\theta}=
  \,^{r_j}\!\Hat{R}_{r_\theta}\Hat{R}_{-\alpha_j}, 
\qquad\text{for $j\neq i$}; \\
\label{vanish_2q}
  \Hat{R}_{r_\theta}=\Hat{R}_{\alpha_i} \mathscr{R},
\end{gather}
where $\mathscr{R}$ is a product of some $R$-matrices.
\end{lm}
\begin{proof}
We have
$r_j\,t_{-\nu(\theta^\vee)}=t_{-\nu(\theta^\vee)}\,r_j$ and
$r_j\,r_0=r_0\,r_j$ for $j\neq i$,
since $\alpha_i$ is the unique vertex connected to $\alpha_0$.
Then $r_j$ and $r_\theta=t_{-\nu(\theta^\vee)}\,r_0$ commute,
which implies $\ell(r_j\,r_\theta)=\ell(r_\theta)+1$ and thus
\eqref{vanish_1q}. 
The form \eqref{vanish_2q} follows from the fact that
$\ell(r_i\,t_{-\nu(\theta^\vee)})=\ell(t_{-\nu(\theta^\vee)})-1$
 implies 
$\ell(r_i\,r_\theta)=\ell(r_\theta)-1$ and
the exchange condition.
\end{proof}
  
Let $\overset{\circ}{W}_i$ be the parabolic subgroup generated 
by $\{r_j | j\in\overset{\circ}{I}, j\neq i\}$ and
$\mathcal{V}_i$ the $\overset{\circ}{W}_i$-invariant
subspace of $\mathcal{M}$.

\begin{lm}
The operator $\Hat{R}_{r_\theta}\,t_{-\nu(\theta^\vee)}$ 
maps $\mathcal{V}_i$ to $\mathcal{V}$
and
the operator
$\Hat{R}_{-\alpha_0}$, 
$\mathcal{V}$ to $\mathcal{V}_i$.
\end{lm}
\begin{proof}
The former statement can be shown in the same way  
as Theorem \ref{symmetric},
and the latter, directly. 
\end{proof}
\begin{thm}
\begin{multline}
\label{qmin}
\Hat{Y}^{-\nu(\theta^\vee)}|_\mathcal{V}= \\
\frac{1}{\bigl|\overset{\circ}{W}_{\nu(\theta^\vee)}\bigr|}
\sum_{w\in\overset{\circ}{W}}w\Biggl(
\Bigl(
\prod_
{\begin{subarray}{c}
\alpha\in\overset{\circ}{\Delta}_+ \\ 
\langle\alpha,\theta^\vee\rangle>0 
\end{subarray}}
H_\alpha(\mu_\alpha)\,
\Bigr)
\Bigl(
H_{a_0^{-1}(\theta+\delta)}(\mu_{\alpha_0})\,t_{-\nu(\theta^\vee)}-
H_{a_0^{-1}(\theta+\delta)}((\overset{\circ}{\rho}_\mu|\theta))
\Bigr)
\Biggr)\Biggr|_\mathcal{V}.
\end{multline}
\end{thm}
\begin{proof}
The explicit form of $\Hat{R}_{r_\theta}\,t_{-\nu(\theta^\vee)}$ 
on $\mathcal{V}_i$
can be computed in a similar way
to Theorem \ref{form_min}.
Since 
$Y^{-\nu(\theta^\vee)}=
R_{r_\theta}\,t_{-\nu(\theta^\vee)}\,R_{-\alpha_0}$,
we obtain the form \eqref{qmin}.
\end{proof}
The operator \eqref{A2_2l} is actually \eqref{qmin} of type
$A^{(2)}_{2l}$, where
the terms without translations are gathered by use of identities of
the theta functions.
In \cite{KomHik97a,KomHik97b},
we calculated the explicit forms of
$\Hat{Y}^{-\lambda_j}|_\mathcal{V}$ for all $j\in\overset{\circ}{I}$
in this root system.
The operator \eqref{qmin}
in the affine root systems of type
$E^{(1)}_8$, $F^{(1)}_4$, $G^{(1)}_2$ and $A^{(2)}_{2l}$
should be compared to 
the Macdonald(-Koornwinder) operator $D_{\theta^\vee}$
of type $E_8$, $F_4$, $G_2$ and $BC_l$ respectively, 
while the operator \eqref{minuscule_op} 
in the rest root systems of type $X^{(1)}_l$
to $E_{\nu^{-1}(\lambda_i)}$
of type $X_l$ \cite{Mac88,Koorn92}.

In order to investigate a general $\Hat{Y}^\lambda$,
let us
define a partial order in $\widehat{M}_-$.
We remark that this partial order 
is different from that in
the affine Hecke algebras.
\begin{defn}
  Let $\lambda,\lambda'\in\widehat{M}_-$.
We write $\lambda\succeq\lambda'$ if $\ell(t_\lambda)>\ell(t_{\lambda'})$ or
$\lambda=\lambda'$.
\end{defn}
For an arbitrary weight $\lambda\in\widehat{M}_-$,
we have the ``leading term'' of $\Hat{Y}^{\lambda}$
with respective to the order $\succ$.
\begin{thm}
 Let $\lambda\in\widehat{M}_-$. Then we have
  \begin{equation}
  \label{expansion_Y}
  \Hat{Y}^\lambda|_\mathcal{V}=
  \frac 1{\bigl|\overset{\circ}{W}_\lambda\bigr|}
  \sum_{w\in\overset{\circ}{W}} w 
  \Bigl(g^\lambda_\lambda t_{\lambda}+
  \sum_{\lambda\succ\lambda'}g^\lambda_{\lambda'}t_{\lambda'}
  \Bigr) 
  \Bigr|_\mathcal{V},
  \end{equation}
where 
$g^\lambda_{\lambda'}\in\mathcal{M}$.
Especially we have 
$g^\lambda_\lambda=\prod_{\alpha\in\Delta_{t_\lambda}} H_\alpha(\mu_\alpha)$.
\end{thm}
\begin{proof}
Because $\Hat{Y}^{\lambda}$ is $\overset{\circ}{W}$-invariant,
it is sufficient to calculate the coefficients of the translations 
of antidominant weights.
A translation $t_{\lambda'}$, $\lambda'\in\widehat{M}_-$
in the expansion of $\Hat{Y}^{\lambda}$
appears as $w t_{\lambda}=t_{\lambda'} \overset{\circ}{w}$
where $\overset{\circ}{w}\in\overset{\circ}{W}$ and 
\begin{gather}
w=r_{\alpha\{p\}}\ldots r_{\alpha\{1\}}\in W,
\\
\alpha\{q\}=\alpha^{m_q}, \qquad
1\leq m_p < m_{p-1} < \ldots < m_2 < m_1 \leq \ell(t_{\lambda}).   
\end{gather}
Then 
$\ell(t_\lambda)\geq \ell(w t_\lambda)=\ell(t_{\lambda'} \overset{\circ}{w} )=
\ell(t_{\lambda'})+\ell(\overset{\circ}{w})$, which implies
$\ell(t_\lambda)>\ell(t_{\lambda'})$ if $w\neq id$.
Hence the expression \eqref{expansion_Y}.
\end{proof}
\begin{thm}[cf.~\cite{Chere92a}]
  $\{\Hat{Y}^{-\lambda_i}|i\in\overset{\circ}{I}\}$ are algebraically independent on $\mathcal{V}$.
\end{thm}
\begin{proof}
Consider $Y=\sum_{\lambda} a_\lambda Y^{\lambda}\in\mathcal{S}$ with
$a_\lambda\in\mathbb{C}$. 
Let $M_Y$ be the set of all the maximal 
antidominant weights 
in the expansion of $\Hat{Y}$ on $\mathcal{V}$.
Let $M'_Y:=\cup_{\lambda\in M_Y}
\{\lambda'\in\widehat{M}_-|\lambda'\preceq\lambda\}$.
Then we have
\begin{equation}
  \Hat{Y}|_\mathcal{V}=\sum_{\lambda\in M_Y} 
  \sum_{w\in\overset{\circ}{W}}w 
  \bigl(
  a_\lambda g^\lambda_\lambda t_\lambda 
+ \text{lower terms ($\lambda'\prec\lambda$)}
  \bigr)\bigr|_\mathcal{V}
\end{equation}
Fix $\lambda\in M_Y$.
There exists $h_0\in\mathfrak{h}$ such that 
\begin{equation}
  \{wh_0-h_0 | w\in\overset{\circ}{W}\}\cap
  \{\kappa^{-1}\langle\delta,h_0\rangle 
      (\nu^{-1}(w \lambda')-\nu^{-1}(w' \lambda'))
      | w,w'\in\overset{\circ}{W}, \lambda'\in M'_Y\}=\emptyset
\end{equation}
and $g^\lambda_\lambda(h_0)\neq 0$.
Suppose $\Hat{Y}f=0$ for all $f\in\mathcal{V}$. Then $(\Hat{Y}f)(h_0)=0$.
Since
$\{(t_{w\lambda'}f)(h_0) |
\lambda'\in M'_Y, w\in\overset{\circ}{W}\}$
 can be made arbitrary for
suitable $f\in\mathcal{V}$,
it follows that
$a_\lambda=0$ and
hence the result.
\end{proof}
\begin{cor}
$\mathcal{S}\simeq\mathbb{C}[T_{\widehat{M}_-}]$.   
\end{cor}

\section{Action on Theta Functions of Level $k$}
\label{sec:theta}
The aim of this section is 
to show that the operators $\Hat{Y}^\lambda$ in the previous sections
act on $(\widetilde{Th}{}^k)^{\overset{\circ}{W}}$,
the $\overset{\circ}{W}$-invariant space 
of the theta functions of level $k$ or the space of the characters. 
To be more precise, 
we identify $\Hat{Y}^\lambda$
with an operator on $(\widetilde{Th}{}^k)^{\overset{\circ}{W}}$
by restricting the domain. 
We regard this space both as a $\mathbb{C}$-vector space and
as an $\mathscr{O}$-module.
The basic idea is from \cite{FelPas94,KomHik96a},
where the matrix elements of
Belavin's $\mathbb{Z}_k$-symmetric elliptic $R$-matrix 
and associated $K$-matrices
are calculated.
Now it is turned out that
they treat
the elliptic difference operators of
type $A^{(1)}_1$ or $A^{(2)}_2$.

First
let us outline our strategy.
Since the representation $\pi$ in Theorem \ref{Thm_R_2} does not
preserve $\widetilde{Th}{}^k$ 
for general $\xi\in\overset{\circ}{\mathfrak{h}}{}^*$,
we introduce another representation $\Bar{\pi}$
which always preserve this space.
The images of $\mathcal{S}$ by $\pi$ and $\Bar{\pi}$
coincide when we set 
$\xi=-\overset{\circ}{\rho}_\mu$.
As was shown, $\pi(\mathcal{S})$ at this value preserves
$\overset{\circ}{W}$-invariant subspace, so does $\Bar{\pi}(\mathcal{S})$.
On the other hand, $\Bar{\pi}(\mathcal{S})$ preserves
$\widetilde{Th}{}^k$ by construction, so does $\pi(\mathcal{S})$.
Therefore we can deduce that $\pi(\mathcal{S})=\Bar{\pi}(\mathcal{S})$
acts on $(\widetilde{Th}{}^k)^{\overset{\circ}{W}}$.

Let $h^\vee_\mu:=(\overset{\circ}{\rho}_\mu|\theta)+\mu_{\alpha_0}
=\sum_{i\in I}\mu_{\alpha_i}a_i^\vee$
and 
$\Xi:=\frac{\xi+\overset{\circ}{\rho}_\mu}{h^\vee_\mu}$.
Throughout this section, we fix $\kappa=\frac{h^\vee_\mu}{k}$
though some of the following statements do not require this condition.

We extend the action of $t_\lambda$ on $\mathcal{M}$
for arbitrary $\lambda\in\overset{\circ}{\mathfrak{h}}{}^*$ by
$(t_\lambda f)(h):=f(t^\iota_{-\lambda} h)$. 
Let $\Bar{R}_\alpha\in\End_\mathbb{C}\mathcal{M}$ be defined by
\begin{equation}
\Bar{R}_\alpha:=t_{\epsilon^1_\alpha} 
  \Hat{R}_{\overline{\alpha}}\,t_{\epsilon^2_\alpha},
\end{equation}
where 
$\epsilon^1_\alpha
 :=\frac{1}{h^\vee_\mu}(-\frac 12\mu_\alpha\overline{\alpha}-\xi+\eta_\alpha)$,
$\epsilon^2_\alpha
 :=\frac{1}{h^\vee_\mu}(-\frac 12\mu_\alpha\overline{\alpha}+\xi-\eta_\alpha)$,
and $\eta_\alpha\in\overset{\circ}{\mathfrak{h}}{}^*$ is taken
arbitrary such that
$\langle\eta_\alpha,\overline{\alpha^\vee}\rangle=0$.
Then
$\Bar{R}_\alpha$ does not depend on the choice of $\eta_\alpha$ and
thus is well-defined.
According to our plan,
we show that this operator acts on $\widetilde{Th}{}^k$.
\begin{prop}
\label{prop:R_end_Th}
For arbitrary $\xi\in\overset{\circ}{\mathfrak{h}}{}^*$,
  $\Bar{R}_\alpha\in\End_{\mathscr{O}}(\widetilde{Th}{}^k)$. 
\end{prop}
\begin{proof}
We note that 
\begin{equation}
  t^\iota_\alpha t_\beta = t_\beta t^\iota_\alpha f^{(\alpha|\beta)},
\end{equation}
where $\alpha\in\widehat{M}$, $\beta\in M$ and
$f(\lambda):=\lambda-\kappa\langle\lambda,K\rangle\delta$ \cite{Chere95b}.
By using this relation and the condition \eqref{cond_theta},
we can check the behavior under the action of the Heisenberg group 
(see the Appendix) and the holomorphy on the domain $Y$. 
Then we see
 $\Bar{R}_\alpha\in\End_{\mathbb{C}}(\widetilde{Th}{}^k)$.
Since $\widehat{W}$ fixes $\delta$,
we have the proof.
\end{proof}

Here we shall make crucial steps to the main statement.
\begin{lm}
\label{lm:eta_j}
Let $\Hat{w}=r_{i_1}\ldots r_{i_\ell}\omega\in\widehat{W}$ 
be a reduced expression.
Let
\begin{equation}
\eta_n:=-\overset{\circ}{\rho}_\mu+
    \sum_{m=1}^{n-1}\nu_m\overline{\alpha^m}+
    \frac 12 \nu_n\overline{\alpha^n},    
\end{equation}
where 
$\Delta_{\Hat{w}}=\{
\alpha^1=\alpha_{i_1},\alpha^2=r_{i_1}(\alpha_{i_2}),\ldots,
\alpha^{\ell}=w r_{i_\ell}(\alpha_{i_\ell})\}$,
\begin{equation}
  \nu_n:=
  \begin{cases}
    \mu_n,&\text{if $\alpha_{i_n}\neq\alpha_0$}; \\
    -(\overset{\circ}{\rho}_\mu|\theta),\qquad
   &\text{if $\alpha_{i_n}=\alpha_0$},
  \end{cases}
\end{equation}
and $\mu_n:=\mu_{\alpha^n}$.
Then 
$\langle\eta_n,\overline{(\alpha^{n})^\vee}\rangle=0$.
\end{lm}
\begin{proof}
First observe that if
$\alpha^n=r_{i_1}\cdots r_{i_{n-1}}\alpha_{i_n}$,
then 
$\overline{\alpha^n}=
\Bar{r}_{i_1}\cdots \Bar{r}_{i_{n-1}}
\overline{\alpha_{i_n}}$ and 
$\overline{(\alpha^n)^\vee}=
\Bar{r}_{i_1}\cdots \Bar{r}_{i_{n-1}}
\overline{\alpha_{i_n}^\vee}$,
where
$\Bar{r}_i:=r_i$ for $i\neq 0$ and 
$\Bar{r}_0:=r_\theta$.
\begin{align}
\langle-\overset{\circ}{\rho}_\mu,\overline{(\alpha^n)^\vee}\rangle
&= \langle-\overset{\circ}{\rho}_\mu,
    \Bar{r}_{i_1} \cdots \Bar{r}_{i_{n-1}}
     \overline{\alpha_{i_n}^\vee}\rangle \notag \\
&= \langle-\overset{\circ}{\rho}_\mu+\nu_1\overline{\alpha^1},
    \Bar{r}_{i_2} \cdots \Bar{r}_{i_{n-1}}
    \overline{\alpha_{i_n}^\vee}\rangle \notag \\
&= \langle-\overset{\circ}{\rho}_\mu,
    \Bar{r}_{i_2} \cdots \Bar{r}_{i_{n-1}}
    \overline{\alpha_{i_n}^\vee}\rangle- 
    \nu_1\langle\overline{\alpha^1},\overline{(\alpha^n)^\vee}\rangle \notag \\
&\ \vdots \notag \\
&= - \sum_{m=1}^{n-1} 
\nu_m\langle\overline{\alpha^m},\overline{(\alpha^n)^\vee}\rangle-\nu_n. 
\end{align}
Then we have
\begin{equation}
\langle-\overset{\circ}{\rho}_\mu+
    \sum_{m=1}^{n-1}\nu_m\overline{\alpha^m}+
    \frac 12 \nu_n\overline{\alpha^n},\overline{(\alpha^n)^\vee}\rangle=0.  
\end{equation}
\end{proof}
\begin{prop}
\label{prop:Bar_eq_Hat}
Let $\Hat{w}=r_{i_1}\cdots r_{i_\ell}\omega\in\widehat{W}$ be a reduced expression. Then
\begin{equation}
\label{Bar_eq_Hat}
\Bar{R}_{\alpha^1}\Bar{R}_{\alpha^2}\cdots
\Bar{R}_{\alpha^\ell}=
t_{-\Xi}
\Hat{R}_{\alpha^1}\Hat{R}_{\alpha^2}\cdots
\Hat{R}_{\alpha^\ell}t_\lambda t_{\Xi},
\end{equation}
where 
$\lambda
=-\frac{1}{h^\vee_\mu}
\sum_{n=1}^{\ell}\mu_n\overline{\alpha^n}
=-\frac{1}{h^\vee_\mu}
\sum_{\alpha\in\Delta_{\Hat{w}}}\mu_{\alpha}\overline{\alpha}$.
\end{prop}
\begin{proof}
We set $\eta_{\alpha^n}=\eta_n$ obtained in Lemma \ref{lm:eta_j}
and set
\begin{gather}
\epsilon_0
:=\epsilon^1_{\alpha^1}
=-\Xi+\frac 1{2h^\vee_\mu}(\nu_1-\mu_1)\overline{\alpha^1},
\\
\epsilon_n
:=\epsilon^2_{\alpha^n}+\epsilon^1_{\alpha^{n+1}}
=\frac{1}{h^\vee_\mu}
 \Bigl(\frac 12(\nu_n-\mu_n)\overline{\alpha^n}+
 \frac 12(\nu_{n+1}-\mu_{n+1})\overline{\alpha^{n+1}}\Bigr),
 \quad 1\leq n\leq \ell-1,
\\
\epsilon_\ell
:=\epsilon^2_{\alpha^\ell}
=\Xi+\frac 1{h^\vee_\mu}\Bigl(
 -\sum_{m=1}^{\ell-1}\nu_m\overline{\alpha^m}-
    \frac 12 (\nu_\ell+\mu_\ell)\overline{\alpha^\ell}
 \Bigr).
\end{gather}
Because 
$r_{i_1}\cdots r_{i_\ell}\omega$ is a reduced expression,
we have for $1\leq n\leq \ell-1$
\begin{equation}
  \epsilon_n=
  \begin{cases}
    -\frac 12\overline{\alpha^n},\qquad&\text{if }\alpha_{i_n}=\alpha_0;\\
    -\frac 12\overline{\alpha^{n+1}},&\text{if }\alpha_{i_{n+1}}=\alpha_0;\\
    0,&\text{otherwise}.\\
  \end{cases}
\end{equation}
Let $\Bar{w}_n=\Bar{r}_{i_1}\cdots\Bar{r}_{i_n}$.
If $\alpha_{i_n}=\alpha_0$, then
$
t_{(-\overline{\alpha^n}/2)}\Hat{R}_{(\overline{\alpha^n})}
t_{(-\overline{\alpha^n}/2)}
=\Hat{R}_{(\Bar{w}_{n-1}\alpha_0)}
t_{(\Bar{w}_{n-1}\nu(\theta^\vee))}
$ and
if $\alpha_{i_n}\neq\alpha_0$, then
$\Hat{R}_{(\overline{\alpha^n})}
=\Hat{R}_{(\Bar{w}_{n-1}\alpha_{i_n})}$.
By using the identity
\begin{equation}
\alpha^n=r_{i_1}\cdots r_{i_{n-1}}\alpha_{i_n}=
\Bigl(
\prod_{\substack{m<n \\ \alpha_{i_m}=\alpha_0}}
t_{(\Bar{w}_{m-1}\nu(\theta^\vee))}
\Bigr)
\Bar{w}_{n-1}\alpha_{i_n},
\end{equation}
we arrive at \eqref{Bar_eq_Hat}.
\end{proof}
Apply this proposition to
an element that has two reduced expressions of the form
\begin{equation}
\Hat{w}=\underbrace{r_i r_j r_i \ldots}_{m_{ij}\text{ factors}} 
 = \underbrace{r_j r_i r_j \ldots}_{m_{ij}\text{ factors}}, 
\end{equation}
for $i\neq j\in I$.
Then the relation 
\begin{equation}
\Bar{R}_{\alpha_i} \Bar{R}_{r_i \alpha_j} \Bar{R}_{r_i r_j \alpha_i}\cdots
=
\Bar{R}_{\alpha_j} \Bar{R}_{r_j \alpha_i} \Bar{R}_{r_j r_i \alpha_j}\cdots
\end{equation}
immediately follows. 
Regarding
$w\Pi$ for $w\in\overset{\circ}{W}$, 
as a set of fundamental roots in
Lemma \ref{lm:eta_j} and 
Proposition \ref{prop:Bar_eq_Hat},
we have proved the following theorem:
\begin{thm}
\label{Thm_R_3}
The map $\Bar{\pi}:R_\alpha \mapsto \Bar{R}_\alpha, t_\lambda \mapsto \Id_\mathcal{M}$
induces a homomorphism from $\mathcal{R}$ to $\End_\mathbb{C}\mathcal{M}$
and $\End_{\mathscr{O}}(\widetilde{Th}{}^k)$.
\end{thm}

For
$\lambda\in\widehat{M}_-$,
we set $\Bar{Y}^{\lambda}:=\Bar{\pi}(Y^{\lambda})=
\Bar{R}_{\alpha^1}\Bar{R}_{\alpha^2}\cdots
\Bar{R}_{\alpha^\ell}\in\End_{\mathscr{O}}(\widetilde{Th}{}^k)$. 
Now we are in position to prove 
the main theorem fully stated as follows:
\begin{thm}
\label{thm:main}
Let $\kappa=\frac{h^\vee_\mu}{k}$ and
$\xi=-\overset{\circ}{\rho}_\mu$. Then
$\Hat{Y}^{\lambda}=\Bar{Y}^{\lambda}\in
\End_{\mathscr{O}}
\bigl((\widetilde{Th}{}^k)^{\overset{\circ}{W}}\bigr)$.
\end{thm} 

By Proposition \ref{prop:Bar_eq_Hat}, we have already shown
\begin{equation}
 \Bar{Y}^\lambda
=t_{-\Xi} \Hat{R}_{\alpha^1}\cdots \Hat{R}_{\alpha^\ell} 
 t_{\lambda'} t_{\Xi}
=t_{-\Xi} \Hat{R}_{t_\lambda}
 t_{\lambda'} t_{\Xi}
  ,
\end{equation}
where $\lambda'=-\frac{1}{h^\vee_\mu}
\sum_{n=1}^{\ell}\mu_n\overline{\alpha^n}
=-\frac{1}{h^\vee_\mu}
\sum_{\alpha\in\Delta_{t_{-\lambda}}}\mu_{\alpha}\overline{\alpha}$.
Since $\Xi=0$ if we set $\xi=-\overset{\circ}{\rho}_\mu$,
we have only 
to show that $\lambda'=\lambda$.

Due to the formulas
\eqref{property}, 
we have another description of $\lambda'$
which can be regarded as an image of $\lambda$ by some linear map:
\begin{equation}
\label{sum_of_roots}
-\sum_{\alpha\in\Delta_{t_{-\lambda}}}\mu_{\alpha}\overline{\alpha}
=
  \begin{cases}
   \displaystyle{
   \sum_{\alpha\in\overset{\circ}{\Delta}_+}\frac{1}{\gamma_\alpha}
   \mu_\alpha(\alpha|\lambda)\alpha
   },
  &\text{if $A$ is not of type $A^{(2)}_{2l}$};\\[2mm]   
   \displaystyle{
   \sum_{\alpha\in\overset{\circ}{\Delta}_+}\frac{1}{\gamma_\alpha}
   \mu_\alpha(\alpha|\lambda)\alpha +
   \frac 14 \mu_{\alpha_0}
   \sum_{\alpha\in(\overset{\circ}{\Delta}_+)_l}
   (\alpha|\lambda)\alpha
   },
  \quad
  &\text{if $A$ is of type $A^{(2)}_{2l}$}.
  \end{cases}
\end{equation}

\begin{lm}
\label{lm:Schur}
Let 
$L:\overset{\circ}{\mathfrak{h}}{}^*\rightarrow
\overset{\circ}{\mathfrak{h}}{}^*$ be a linear map defined by
$L:\lambda\mapsto\frac 12\sum_{\alpha\in\overset{\circ}{\Delta}}
\nu_\alpha(\alpha|\lambda)\alpha$ where $\nu_\alpha$ is 
$\overset{\circ}{W}$-invariant constant.
Then $L=a\Id_{\overset{\circ}{\mathfrak{h}}{}^*}$ where $a\in\mathbb{C}$.
\end{lm}
\begin{proof}
We see $L\in\End_{\mathbb{C}[\overset{\circ}{W}]}
(\overset{\circ}{\mathfrak{h}}{}^*)$.
Since $\mathbb{C}[\overset{\circ}{W}]$ acts 
on $\overset{\circ}{\mathfrak{h}}{}^*$ irreducibly,
the statement follows from Schur's lemma.
\end{proof}
By this lemma, we see that
$-\sum_{\alpha\in\Delta_{t_{-\lambda}}}\mu_{\alpha}\overline{\alpha}
=a\,\lambda$
for some $a\in\mathbb{C}$.
The following proposition 
completes the proof of Theorem \ref{thm:main}.  
\begin{prop}
  $-\sum_{\alpha\in\Delta_{t_{-\lambda}}}\mu_{\alpha}\overline{\alpha}=
  h^\vee_\mu\,\lambda$.
\end{prop}
\begin{proof}
Let $L$ be a linear map defined in the right hand side 
of \eqref{sum_of_roots}.
Owing to Lemma \ref{lm:Schur}, 
we can evaluate the factor $a$
at any element of $\overset{\circ}{\mathfrak{h}}{}^*$. 
Recall that every root system has a quasi-minuscule weight
$\nu(\theta^\vee)$, 
whose properties 
we have already investigated.
\begin{itemize}
\item $A$ is not of type $A^{(2)}_{2l}$ 
  \begin{equation}
L(\nu(\theta^\vee))
=\sum_{\alpha\in\overset{\circ}{\Delta}_+}
  \frac{1}{\gamma_\alpha}
  \mu_\alpha\langle\alpha,\theta^\vee\rangle\alpha
=\sum_{\substack{
  \alpha\in\overset{\circ}{\Delta}_+ \\
  \langle\alpha,\theta^\vee\rangle\neq 0
  }}
  \mu_\alpha \alpha+\mu_\theta \theta=a\,\nu(\theta^\vee) ,
  \end{equation}
where we have used
Lemma \ref{property_theta}.
By applying $(.|\theta)$ in the last equality,
we obtain
\begin{equation}
  a=\frac 12 \sum_{\alpha\in\overset{\circ}{\Delta}_+}
  \mu_\alpha (\alpha|\theta)+\mu_\theta
 =(\overset{\circ}{\rho}_\mu|\theta)+\mu_{\alpha_0}.
\end{equation}
\item  $A$ is of type $A^{(2)}_{2l}$ \\
In a similar manner, we have
\begin{equation}
L(\nu(\theta^\vee))
=\sum_{\substack{
  \alpha\in\overset{\circ}{\Delta}_+ \\
  \langle\alpha,\theta^\vee\rangle\neq 0
  }}
  \mu_\alpha \alpha +
\sum_{\substack{
  \alpha\in(\overset{\circ}{\Delta}_+)_l \\
  \langle\alpha,\theta^\vee\rangle\neq 0
  }}
  \frac 12 \mu_{\alpha_0} \alpha =a\, \nu(\theta^\vee),
\end{equation}
and consequently
\begin{equation}
a=(\overset{\circ}{\rho}_\mu|\theta)
 +\frac 14 \mu_{\alpha_0} 
  (\rho_l|\theta)
=(\overset{\circ}{\rho}_\mu|\theta)+\mu_{\alpha_0},
\end{equation}
where $\displaystyle{\rho_l=
\sum_{\alpha\in(\overset{\circ}{\Delta}_+)_l}\alpha=
2\overline{\Lambda_l}}$.
\end{itemize}

In any case,
$L(\nu(\theta^\vee))=h^\vee_\mu \,\nu(\theta^\vee)$ and we have
$-\sum_{\alpha\in\Delta_{t_{-\lambda}}}\mu_{\alpha}\overline{\alpha}=
L(\lambda)=h^\vee_\mu\,\lambda$, as required.
\end{proof}

Note that we also showed that
\begin{equation}
  \sum_{\alpha\in\overset{\circ}{\Delta}}(\lambda|\alpha)(\mu|\alpha)
  =2h^\vee(\lambda|\mu),\quad
\text{for $\lambda,\mu\in\overset{\circ}{\mathfrak{h}}{}^*$},
\end{equation}
in the nontwisted root systems.
See Corollary 8.7 of \cite{Kac}.

\section{Concluding Remarks}

We constructed mutually commuting difference operators
 by means of the root algebras.
Since the operator is represented 
in a single product of affine $R$-matrices,
we had only to pursue the image of each $R$-matrix
and 
therefore suceeded in proving that 
they act on the characters of the 
irreducible representations
of the affine Lie algebras.
However, the procedure of the diagonalization has yet to be solved.
Prior to this difficult problem,
we may need 
to show the selfadjointness on the space of the characters 
with respect to some inner product,
since there is no certainty that they can be diagonalized.
In $A^{(1)}_2$ case, the selfadjointness was established in \cite{Skl83}
for arbitrary level of positive integer. 

Since this operator was originally introduced
as a quantum many-body system,
the selfadjointness should be also an important problem
in this sense.
In the trigonometric case, we readily see that
the Macdonald operators are essentially selfadjoint 
on the polynomials of exponential since the operators
are diagonalized in terms of the Macdonald polynomials.
In the elliptic case, however, 
this problem is less investigated.
See, for example,
\cite{Ruij99a,Ruij99b} where the two-body system is extensively studied
by constructing the explicit eigenvectors,
or \cite{Kom98a} where
the extensibility to 
positive selfadjoint operators is
shown by introducing a certain measure on a torus.
These systems correspond to negative levels 
in terms of the affine Lie algebras
and if we treat positive level cases, 
the measure includes discrete parts.

We hope the construction developed in this paper 
shed light on these problems.


\section*{Acknowledgment}
The author expresses his sincere gratitude to 
Prof.~Masaki Kashiwara and Prof.~Tetsuji Miwa 
who kindly allowed him to speak about his research.
Also thanks are due to 
Prof.~Atsuo Kuniba,
Prof.~Jun'ichi Shiraishi, Prof.~Junji Suzuki, Prof.~Yuji Yamada,
Dr.~T. H. Baker, Dr.~Goro Hatayama, Prof.~Koji Hasegawa, Dr.~Takeshi Ikeda,
Dr.~Tetsuya Kikuchi,
Dr.~Kazuhiro Hikami 
and Dr.~Akinori Nishino
for fruitful discussions and helpful comments.
He would like to thank Prof.~Miki Wadati 
for kind interests in this work.
He is a Research Fellow of the Japan
Society for the Promotion of Science.

\appendix
\section{Fundamental Functions and Identities}
We define an action of 
$n=(v,\lambda,u)\in\overset{\circ}{\mathfrak{h}}\times
\overset{\circ}{\mathfrak{h}}{}^*\times \mathbb{C}$ 
on a holomorphic function $F$
on $Y$ by
\begin{equation}
  (n F)(h):=F(t_{-\lambda}(h)-2\pi i v - 
        (u+\pi i\langle\lambda,v\rangle)K).
\end{equation}
\begin{defn}
The Heisenberg group is $N_{\mathbb{Z}}=\{(v,\lambda,u)
  \in\overset{\circ}{\mathfrak{h}}\times
  \overset{\circ}{\mathfrak{h}}{}^*\times i\mathbb{R} |
  v\in\overset{\circ}{Q}{}^\vee, \lambda\in M,
  u+\pi i \langle\lambda,v\rangle\in 2\pi i\mathbb{Z}\}$ with multiplication:
\begin{equation}
(v,\lambda,u)(v',\lambda',u')
:=
(v+v',\lambda+\lambda',u+u'+     
\pi i (\langle\lambda',v\rangle-\langle\lambda,v'\rangle) ).  
\end{equation}
\end{defn}
\begin{defn}
Fix a nonnegative integer $k$. A theta function of level $k$ is 
a holomorphic function $F$ on the domain $Y$ such that the following
 two conditions hold:
 \begin{gather}
n(F)=F \qquad \text{for all $n\in N_\mathbb{Z}$}; \\
n(F)=e^{-ka} F \qquad \text{for all $n=(0,0,a) \in (0,0,\mathbb{C})$}.
\end{gather}
\end{defn}
Let $\widetilde{Th}{}^k$ denote the vector space over $\mathbb{C}$
of the theta functions of level $k$.
It is known that $\mathscr{O}:=\widetilde{Th}{}^0$ is the set of
holomorphic functions of $\langle\delta,h\rangle$.

For $\lambda\in\mathfrak{h}^*$ such that level$(\lambda)=k>0$, we set
\begin{equation}
  \Theta_\lambda:=e^{-\frac{|\lambda|^2}{2k}\delta}
\sum_{t\in T_M} e^{t(\lambda)}.
\end{equation}
It is known that $\{\Theta_\lambda|\text{ level$(\lambda)=k$}\}$
is an $\mathscr{O}$-basis of $\widetilde{Th}{}^k$. 

Consider the root system of type $A^{(1)}_1$.
Then $\Pi=\{\alpha_0,\alpha_1\}$;
$M=\mathbb{Z}\alpha_1$;
$(\alpha_1|\alpha_1)=2$.
We have 4 theta functions of level $2$ for $k\in\mathbb{Z}/4\mathbb{Z}$;
\begin{equation}
\label{theta_4}
    \Theta_{2\Lambda_0+k\overline{\Lambda_1}}
    =
    e(2\Lambda_0)
    \sum_{n \in \mathbb{Z} }
    e
    \Bigl(
      - \frac 12 (2n+\frac k2)^2 \delta
      + (2n+\frac k2) \alpha_1
    \Bigr), 
\end{equation}
where $e(\lambda)(h):=\exp(\langle\lambda,h\rangle)$
for $\lambda\in\mathfrak{h}^*$.
We see that $\text{level}(\rho)=h^\vee=2$.
\begin{equation}
\label{A_rho_2}
  A_\rho=\sum_{w\in\overset{\circ}{W}}\varepsilon(w)\Theta_{w(\rho)} 
        =\Theta_{2\Lambda_0+\frac 12\alpha_1}-\Theta_{2\Lambda_0-\frac 12\alpha_1} 
=e(\rho-\frac 18\delta)\prod_{\alpha\in\Delta_+}(1-e(-\alpha)).
\end{equation}
Motivated by these equations,
we define theta functions
for $\lambda\in\mathfrak{h}^*$ and $\gamma>0$ by
\begin{subequations}
  \begin{gather}
    \vartheta^1(\lambda;\gamma)
     :=
    \sum_{n \in \mathbb{Z} }
    (-1)^n e
    \Bigl(
      - \frac 12 (n+\frac 12)^2 \gamma\delta
      + (n+\frac 12) \lambda
    \Bigr) ,\\
    \vartheta^2(\lambda;\gamma)
    :=
    \sum_{n \in \mathbb{Z} }
    e
    \Bigl(
      - \frac 12 (n+\frac 12)^2 \gamma\delta
      + (n + \frac 12 ) \lambda
    \Bigr) ,\\
    \vartheta^3(\lambda;\gamma)
    :=
    \sum_{n \in \mathbb{Z} }
    e
    \Bigl(
      - \frac 12 n^2 \gamma\delta
      + n \lambda
    \Bigr) ,\\
    \vartheta^0(\lambda;\gamma)
    :=
    \sum_{n \in \mathbb{Z} }
    (-1)^n e
    \Bigl(
      - \frac 12 n^2 \gamma\delta
      + n \lambda
    \Bigr) ,
  \end{gather}
\end{subequations}
and eta function
\begin{equation}
\eta(\delta)
    :=
 e\Bigl(-\frac 1{24}\delta\Bigr)
\prod_{n\in\mathbb{Z}_{\geq 1}}
(1-e(-n\delta)).
\end{equation}
Then we have by \eqref{theta_4} and \eqref{A_rho_2}
\begin{equation}
 \vartheta^1(\lambda;\gamma)=
e\Bigl(
\frac{\lambda}2-\frac 18\gamma\delta\Bigr)
(1-e(-\lambda))
\prod_{\substack{
\lambda'\in\{-\lambda,0,\lambda\} \\
n\in\mathbb{Z}_{\geq 1}
}}
(1-e(-\lambda'-n\gamma\delta)).
\end{equation}
It is well known that
\begin{equation}
  A_\rho=
   e(h^\vee\Lambda_0)f(\delta)
   \prod_{\alpha\in\overset{\circ}{\Delta}_+}
   \vartheta^1(\alpha;\gamma_\alpha),
\end{equation}
for some function $f(\delta)$ which depends only on $\delta$.
We symbolically set 
\begin{equation}
  \vartheta^{1\prime}(0;\gamma):=(\eta(\gamma\delta))^3.
\end{equation}
For $\lambda\in\mathfrak{h}^*$ and $\nu\in\mathbb{C}$, we define
\begin{equation}
  \sigma_\nu(\lambda;\gamma)
  :=
  \frac 
{ \vartheta^1(\lambda-\nu\delta;\gamma) \vartheta^{1\prime} (0;\gamma) }
{ \vartheta^1(\lambda;\gamma) \vartheta^1(-\nu\delta;\gamma) }, \qquad
  \wp^0(\lambda;\gamma)
  :=\Biggl(
  \frac{\vartheta^0(\lambda;\gamma) \vartheta^{1\prime} (0;\gamma) }
  {\vartheta^1(\lambda;\gamma) \, \vartheta^0(0;\gamma) } \Biggr)^2.
\end{equation}

We see that these theta functions are related 
to the classical Jacobi theta functions $\vartheta_j(z;\tau)$,
 Weierstrass function $\wp(z;1,\tau)$ and
 Dedekind eta function $\eta(\tau)$
 as
\begin{gather}
  \vartheta^1(\lambda;\gamma)(h)
  =-i
  \vartheta_1 (\langle \lambda,\overset{\circ}{h} \rangle
       ;\gamma\tau ), \\
  \vartheta^j(\lambda;\gamma)(h)
  =
  \vartheta_j (\langle \lambda,\overset{\circ}{h} \rangle
       ;\gamma\tau ), \\
  \wp^0(\lambda;\gamma)(h)=-\frac{1}{4\pi^2}
  ( \wp (\langle \lambda,\overset{\circ}{h} \rangle;1,\gamma\tau )
   -\wp (\gamma\tau/2;1,\gamma\tau ) ) ,\\
  \eta(\delta)(h)=\eta(\gamma\tau) ,
\end{gather}
where we have set $h=2\pi i ( \overset{\circ}{h} -\tau d + u K )$,
$\overset{\circ}{h}\in\overset{\circ}{\mathfrak{h}}$
and $\tau,u\in\mathbb{C}$.
\bibliographystyle{amsplain}
\bibliography{bibliography}

\end{document}